\documentclass{article}
\usepackage{graphicx,url}
\usepackage{multirow,verbatim}
\usepackage{epsfig,epsf,listings}
\usepackage{amsmath,amssymb,amsfonts,bbm,mathrsfs}
\usepackage{caption}
\usepackage{cite}
\usepackage{float}
\usepackage{mhchem,ulem}
\usepackage{subfigure}
\usepackage{color}
\usepackage{fancyhdr}
\usepackage{epstopdf}

\renewcommand{\xi}{\overline{x}}

\renewcommand{\thispagestyle}[2]{}

\fancypagestyle{plain}{
        \fancyhead{}
        \fancyhead[C]{first page center header}
        \fancyfoot{}
        \fancyfoot[C]{first page center footer}
}
\begin{document}
\thispagestyle{empty}
\setcounter{page}{0}
\begin{Huge}
\begin{center}
Computer Science Technical Report CSTR-{\tt 8} \\
\today
\end{center}
\end{Huge}
\vfil
\begin{huge}
\begin{center}
Azam S. Zavar Moosavi, Adrian Sandu 
\end{center}
\end{huge}

\vfil
\begin{huge}
\begin{it}
\begin{center}
``Approximate Exponential Algorithms to\\ Solve the Chemical Master Equation''
\end{center}
\end{it}
\end{huge}
\vfil

\begin{large}
\begin{center}
Computational Science Laboratory \\
Computer Science Department \\
Virginia Polytechnic Institute and State University \\
Blacksburg, VA 24060 \\
Phone: (540)-231-2193 \\
Fax: (540)-231-6075 \\ 
Email: \url{sandu@cs.vt.edu} \\
Web: \url{http://csl.cs.vt.edu}
\end{center}
\end{large}

\vspace*{1cm}

\begin{tabular}{ccc}
\includegraphics[width=2.5in]{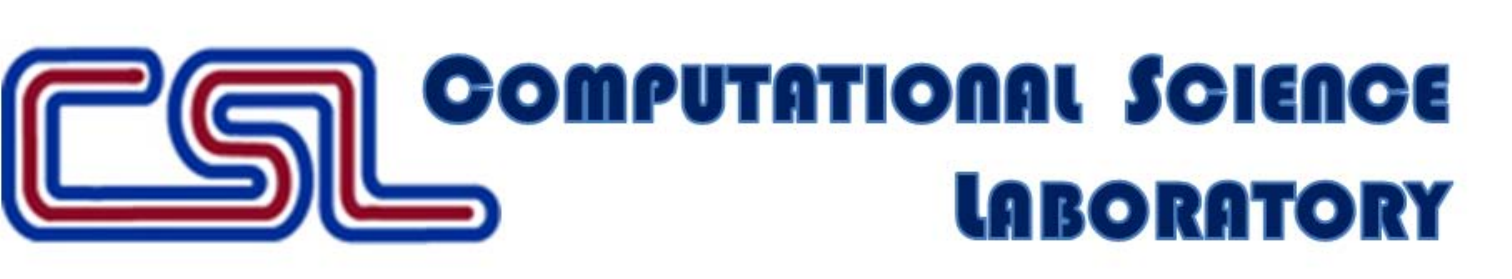}
&\hspace{2.5in}&
\includegraphics[width=2.5in]{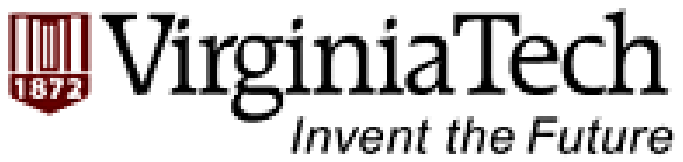} \\
{\bf\em Innovative Computational Solutions} &&\\
\end{tabular}

\newpage

\begin{abstract}
This paper discusses new simulation algorithms for stochastic chemical kinetics that exploit the linearity of the chemical master equation and its matrix exponential exact solution. 
These algorithms make use of various approximations of the matrix exponential to evolve probability densities in time. A sampling of the approximate solutions of the chemical master equation is used to derive accelerated 
stochastic simulation algorithms. Numerical experiments compare the new methods with the established stochastic simulation algorithm and the tau-leaping method.

\paragraph{keywords} Stochastic chemical kinetics, chemical master equation, exact solution, stochastic simulation algorithm, tau-leap.

\end{abstract}

\section{Introduction}
\label{sect:Intro}
In many biological systems the small number of participating molecules make the chemical reactions inherently stochastic. The system state is described by probability densities of the numbers of molecules of different species. The evolution of probabilities in time is described by the chemical master equation (CME) \cite{Gillespie_1977}.
Gillespie proposed the Stochastic Simulation Algorithm (SSA), a Monte Carlo approach that samples from CME \cite{Gillespie_1977}. SSA became the standard method for solving well-stirred chemically reacting systems. However, SSA simulates one reaction and is inefficient for most realistic problems. This motivated the quest for approximate sampling techniques to enhance the efficiency.

The first approximate acceleration technique  is the tau-leaping method \cite{Gillespie_2001} which is able to simulate multiple chemical reactions appearing in a pre-selected time step of length $\tau$. The tau-leap method is accurate if $\tau$ is small enough to satisfy the leap condition,  meaning  that propensity functions remain nearly constant in a time step. The number of firing reactions in a time step is approximated by a Poisson random variable \cite{Kurz_1972_SSA}. 
Explicit tau-leaping method is numerically unstable for stiff systems  \cite{Cao_2004_stability}. Stiffness systems have well-separated ``fast'' and
``slow'' time scales present, and the ``fast modes'' are stable. 
The implicit tau-leap method \cite{Rathinam_2003} overcomes the stability issue  but it has a damping effect
on the computed variances. More accurate variations of the implicit tau-leap method have been proposed  
to alleviate the damping \cite{Gillespie_2003,Gillespie_2001,Sandu_2013_SSA,Cao_2005,Cao_2004,Rathinam_2005}. Simulation efficiency has been increased via parallelization \cite{Sandu_2012_parallel}.

Direct solutions of the CME are computationally important specially in order to estimate moments of the distributions of the chemical species \cite{Burrage_2006_multiscale}. Various approaches to solve the CME are discussed in \cite{Engblom_2006_thesis}.

Sandu has explained the explicit tau-leap method as an exact sampling procedure from an approximate solution of the CME \cite{Sandu_2013_CME}. This paper extends that study and proposes new approximations to the CME solution based on various approximations of matrix exponentials. Accelerated stochastic simulation algorithms are the built by performing exact sampling of these approximate probability densities.

The paper is organized as follows. Section \ref{sect:StochasticChem} reviews the stochastic simulation of chemical kinetics. Section \ref{sect:ApproxExponential} developed the new approximation methods. Numerical experiments to illustrate the proposed schemes are carried out in Section \ref{sect:NumericalExperim}. Conclusions are drawn in Section \ref{sect:Conclusion}.

\section{Simulation of stochastic chemical kinetics}
\label{sect:StochasticChem}
Consider a chemical system in a constant volume container. The system is
well-stirred and in thermal equilibrium at some constant temperature. There
are $N$ different chemical species $S^1,\, \ldots\,,S^N$. Let $X^i(t)$
denote the number of molecules of species $S_i$ at time $t$. The state
vector $x(t)=[X^1(t),\, \ldots\,,X^N(t)]$ defines the numbers of molecules
of each species present at time $t$. The chemical network consists of $M$
reaction channels $R_1,\, \ldots\,,R_M$. Each individual reaction destroys a
number of molecules of reactant species, and produces a number of molecules
of the products. Let $\nu_j^i$ be the change in the number of $S^i$
molecules caused by a single reaction $R_j$. The state change vector $\nu_j
= [\nu_j^{1}, \ldots , \nu_j^{N}]$ describes the change in the entire state
following $R_j$.

A propensity function $a_j(x)$ is associated with each reaction channel $R_j$. 
The probability that one $R_j$ reaction will occur in the next
infinitesimal time interval $[t, t+dt)$ is $a_j(x(t))\cdot dt$. The purpose
of a stochastic chemical simulation is to trace the time evolution of the
system state $x(t)$ given that at the initial time $\bar{t}$ the system is
in the initial state $x\left(\bar{t}\right)$.

\subsection{Chemical Master Equation}
\label{sect:ChemicalMaster}

The Chemical Master Equation (CME) \cite{Gillespie_1977} has complete information about time
evolution of probability of system's state
\begin{equation}\label{eq_CME}
\frac {\partial\mathcal{P}\left(x,t\right)}{\partial t}=\sum_{r=1}^M a_{r}\left(x-v_{r}\right)%
\mathcal{P}\left(x-v_{r},t\right)-a_0\left(x\right)\mathcal{P}\left(x,t\right)\,.
\end{equation}
Let $Q^i$ be the total possible number of molecules of species $S^i$. The
total number of all possible states of the system is: 
\[ \label{eq_Q}
Q=\prod_{i=1}^{N}\left(Q^i+1\right).
\]
We denote by $\mathcal{I}(x)$ the state-space index of state $x=[X^1,\, \ldots\,,X^N]$
\[
\begin{array}{l}
 \mathcal{I}(x) = \left(Q^{N-1}+1\right)\cdots 
\left(Q^1+1\right)\cdot X^N+\cdots  \\ 
 +\left(Q^2+1\right)\left(Q^1+1\right)\cdot X^3+
\left(Q^1+1\right)\cdot X^2+X^1+1 
\end{array}%
\]

One firing of reaction $R_{r}$ changes the state from $x$ to $\bar {x}=x-v_{r}$.
The corresponding change in state space index is: 
\[  \label{eq_d}
\begin{array}{l}
\mathcal{I}(x)-\mathcal{I}\left(x-v_{r}\right)=d_{r},   \\ 
d_{r}=\left(Q^{N-1}+1\right)\cdots\left(Q^1+1\right).v_{r}^N+...\\
\qquad +\left(Q^2+1\right)\left(Q^1+1\right).v_{r}^3+\left(Q^1+1\right).v_{r}^2+v_{r}^1.
\end{array}
\]
The discrete solutions of the CME \eqref{eq_CME} are vectors in the discrete state space, $\mathcal{P}\left(t\right) \in  \mathbb{R}^{Q}$.
Consider  the diagonal matrix  $A_{0}\in \mathbb{R}^{Q \times Q} $ and the Toeplitz matrices $A_{1},\cdots,A_{M}\in \mathbb{R}^{Q \times Q} $ \cite{Sandu_2013_CME}
\[  \label{eq_sumofexponentexact}
({A_{0}})_{i,j}=\left\{ 
\begin{array}{rl}
-a_{0}\left(x_j\right) & \mbox{if $i=j$} ,\\ 
0 & \mbox{if $i \not =j$} ,
\end{array}
\right. \,, \quad
({A_{r}})_{i,j}=\left\{ 
\begin{array}{rl}
a_{r}(x_j) & \mbox{if $i-j=d_{r}$}, \\ 
0 & \mbox{if $i-j \not =d_{r}$},
\end{array}
\right.
\]
as well as their sum $A \in \mathbb{R}^{Q \times Q}$ with entries
\begin{equation}  \label{eq_exact}
A = A_{0} + A_{1} + \dots + A_{M}\,, \quad
A_{i,j}=\left\{ 
\begin{array}{rl}
-a_{0}(x_j) & \mbox{if }i=j\,, \\ 
a_{r}(x_j) & \mbox{if }i-j=d_{r},~ r=1,\cdots,M\,, \\ 
0 & \mbox{otherwise} \,,
\end{array}
\right.
\end{equation}

where $x_j$ denotes the unique state with state space index $j=\mathcal{I}(x_j)$.
In fact matrix A is a square $\left(Q \times Q \right)$ matrix which contains all the propensity values for each possible value of all species or let's say all possible states of reaction system. All possible states for a reaction system consists of $N$ species where each specie has at most $Q^{i}$ $i=1,2,...,N$ value. 

The CME \eqref{eq_CME} is a linear ODE on the discrete state space 
\begin{equation}  \label{eq_cme_mat}
\mathcal{P}' = A \cdot \mathcal{P}\,, \quad  \mathcal{P}(\bar{t}) = \delta_{\mathcal{I}(\bar{x})}\,, \quad
t \ge \bar{t}\,,
\end{equation}
where the system is initially in the known state $x(0)=\bar{x}$ and therefore the initial probability distribution 
vector $\mathcal{P}(0) \in  \mathbb{R}^{Q}$ is equal to one at $\mathcal{I}(\bar{x})$ and is zero everywhere else.
The exact solution of the linear ODE \eqref{eq_cme_mat} is follows: 
\begin{equation}  \label{eq_exact2}
\mathcal{P}\left(\bar{t}+T\right)=\exp\left(T\, A\right)\cdot \mathcal{P}\left(\bar{t}\right) 
= \exp\left(T\, \sum_{r=0}^M A_r\right)\cdot  \mathcal{P}\left(\bar{t}\right)\,.
\end{equation}
%

\subsection{Approximation to Chemical Master Equation}
\label{sect:ApproxChemicalMaster}

Although the CME \eqref{eq_CME} fully describes the evolution of probabilities it is difficult
to solve in practice due to large state space. Sandu \cite{Sandu_2013_CME} considers the following approximation of the CME:
\begin{equation} \label{eq_approx_CME}
\frac {\partial\mathcal{P}\left(x,t\right)}{\partial t}=\sum_{r=1}^M a_{r}\left(\bar{x}\right)%
\mathcal{P}\left(x-v_{r},t\right)- a_0\left(\bar{x}\right)\mathcal{P}\left(x,t\right)
\end{equation}
where the arguments of all propensity functions have been changed from $x$
or $x-v_{j} $ to $\bar{x} $. 
In order to obtain an exponential solution to \eqref{eq_approx_CME} in probability space we consider the diagonal matrix  $\bar{A_{0}}\in \mathbb{R}^{Q \times Q} $ and the Toeplitz matrices $\bar{A_{1}},...,\bar{A_{M}}\in \mathbb{R}^{Q \times Q} $ \cite{Sandu_2013_CME}.
$\bar{A_{r}}$ matrices are square $\left(Q \times Q\right)
$ matrices are built upon the current state of system in reaction system which is against $A_{r}$ matrices that contain all possible states of reaction system.

\begin{equation}  \label{eq_sumofexponent}
(\bar{A_{0}})_{i,j}=\left\{ 
\begin{array}{rl}
-a_{0}\left(\bar{x}\right) & \mbox{if $i=j$} ,\\ 
0 & \mbox{if $i \not =j$} ,
\end{array}
\right. \,, \quad
(\bar{A_{r}})_{i,j}=\left\{ 
\begin{array}{rl}
a_{r}(\bar{x}) & \mbox{if $i-j=d_{r}$}, \\ 
0 & \mbox{if $i-j \not =d_{r}$},
\end{array}
\right.
\end{equation}

together with their sum $\bar{A} = \bar{A_{0}} + \dots + \bar{A_{M}}$.
The approximate CME \eqref{eq_approx_CME} can be written as the linear ODE
\[
\label{eq_sumofexponent2}
\mathcal{P}' = \bar{A} \cdot \mathcal{P}\,, \quad  \mathcal{P}(\bar{t}) = \delta_{\mathcal{I}(\bar{x})}\,, \quad
t \ge \bar{t}\,,
\]
and has an exact solution
\begin{equation}  \label{eq_sumofexponent3}
\mathcal{P}\left(\bar{t}+T\right)=\exp\left(T\, \bar{A}\right)\cdot \mathcal{P}\left(\bar{t}\right)
= \exp\left(T\, \sum_{r=0}^M \bar{A}_r\right)\cdot  \mathcal{P}\left(\bar{t}\right)\,.
\end{equation}
%

\subsection{Tau-leaping method}
\label{sect:TauLeap}

In tau-leap method the number of times a reaction
fires is a random variable from a Poisson distribution with parameter $a_{r}\left(
\bar{x}\right)\tau$.  
Since each reaction fires independently, the probability that each reaction $
R_{r} $ fires exactly $k_{r} $ times, $r=1, 2,\cdots, M $, is the product of $
M $ Poisson probabilities. 
\[
\mathcal{P}\left(K_{1}=k_{1},\cdots,K_{M}=k_{M}\right)= \prod_{r=1}^{M}e^{-a_{r}\left(\bar{x}
\right)\tau}\cdot\frac{\left(a_{r}(\bar{x}\tau\right)^{k_{r}}}{K_{r}!}= e^{-a_{0}\left(\bar{x}
\right)\tau}\cdot 
\prod_{r=1}^{M}\frac{\left(a_{r}\left(\bar{x}\tau\right)\right)^{k_{r}}}{K_{r}!}
\]
Then the state vector after these reactions will change as follows: 
\begin{equation}  \label{eq_tauleap}
X\left(\bar{t}+\tau\right)=\bar{x}+\sum_{r=1}^{M}K_{r}v_{r}
\end{equation}
The probability to go from state $\bar{x}$ at $\bar{t}$ to state $x$ at $
\bar{t}+\tau$, $\mathcal{P}\left(X\left(\bar{t}+\tau\right)\right)=x$, is the sum of all possible
firing reactions which is: 
\[  \label{eq_tauleap1}
\mathcal{P}\left(X,\bar{t}+\tau\right)=e^{-a_{0}\left(\bar{x}\right)T}\cdot\Sigma_{k \in \mathcal{K
}\left(x - \xi\right)} ~ \prod_{r=1}^{M}\frac{\left(a_{r}\left(\bar{x}T\right)\right)^{k_{r}}}{K_{r}!}
\]
Equation \eqref{eq_sumofexponent3} can be approximated by product of each matrix
exponential: 
\begin{equation}  \label{eq__productofexponents}
\mathcal{P}\left(\bar{t}+T\right)=\exp\left(T\bar{A_{0}}\right)\cdot \exp\left(T\bar{A_{1}}
\right)\cdots  \exp\left(T\bar{A_{r}}\right) \cdot \mathcal{P}\left(\bar{t}\right).
\end{equation}
 It has been shown in \cite{Sandu_2013_CME} that the probability given by the tau-leaping method is exactly the probability evolved by the approximate solution \eqref{eq__productofexponents}.
\section{Approximations to the exponential solution}
\label{sect:ApproxExponential}

\subsection{Strang splitting}
\label{sect:StrangSplitting}
In order to improve the approximation of the matrix exponential in 
\eqref{eq__productofexponents} we consider the symmetric Strang splitting \cite{Strang_1968}. For $T=n\tau$ Strang splitting applied to an interval of length $\tau$ leads to the approximation
\begin{equation}  \label{eq_strang}
\mathcal{P}\left(\bar{t}+i \tau\right)=e^{\tau/2 \bar{A}_{r}}\cdots e^{\tau/2 \bar{A}_{1}%
}e^{\tau/2 \bar{A}_{0}}
 \cdot e^{\tau/2 \bar{A}_{1}} \cdots e^{\tau/2 \bar{A}_{r}}\cdot P\left(\bar{t} + (i-1)\tau\right)
\end{equation}
where the matrices $\bar{A_{r}}$ are defined in
\eqref{eq_sumofexponent}.

\subsection{Column based splitting}
\label{sect:ColumnBased}
In column based splitting the matrix $A$  \eqref{eq_exact}  is decomposed in a sum of columns
\[  \label{eq_colbased}
A=\sum_{j=1}^Q A_{j}\,, \quad A_{j}=c_{j}e_{j}^T\,.
\]
Each matrix $A_{j}$ has the same $j$-th column as the matrix $A$, and is zero everywhere else.
Here $c_{j} $ is the $j_{th}$ column of matrix A and $e_{j}$ is the
canonical vector which is zero every where except the $j_{th}$ component. The
exponential of $\tau A_{j}$ is: 
\begin{equation}  \label{eq_colbased2}
e^{\tau A_{j}}=\sum_{k \ge 0}
\frac {\tau^k \left(A_{j}\right)^k}{k!}\,.
\end{equation}
Since $e_{j}^Tc_{j}$ is equal to the $j$-th diagonal entry of matrix A: 
\[  \label{eq_colbased4}
e_{j}^T\,c_{j}=-a_{0}\left(x_{j}\right)
\]
the matrix power $A_{j}^k$ reads
\[  \label{eq_colbased3}
A_{j}^k=c_{j}e_{j}^T \, c_{j}e_{j}^T\, \cdots \, c_{j}e_{j}^T
=  \left(-a_{0}\left(x_{j}\right)\right)^{k-1} c_{j}e_{j}^T = \left(-a_{0}\left(x_{j}\right)\right)^{k-1}A_{j}\,.
\]
Consequently the matrix exponential \eqref{eq_colbased2} becomes 
\[  \label{eq_colbased5}
e^{\tau A_{j}}=I+ \sum_{k \geq 1} \frac {\left(-\tau a_{0}\left(x_{j}\right)\right)^{k-1}}{k!}\left(\tau
A_{j}\right) = I+ S_{j}\, \tau A_{j}\,, \quad S_{j}=\sum_{k \geq 1} \frac {\left(-\tau a_{0}\left(x_{j}\right)\right)^{k-1}}{k!}\,.
\]
We have
\[  \label{eq_colbased7}
e^{\tau A}=e^{\tau \sum_{j=1}^Q A_{j}}\approx \prod_{j=1}^Q e^{\tau A_{j}}
\approx\prod_{j=1}^Q \left(I+S_{j}\tau A_{j}\right)
\]
and the approximation to the CME solution reads 
\[  
\mathcal{P}\left(\bar{t}+i \tau\right)\approx \prod_{j=1}^Q \left(I+S_{j}\tau A_{j}\right)\cdot P\left(\bar{t} + (i-1)\tau\right)\,.
\]
%

\subsection{Accelerated tau-leaping}
\label{sect:AccelTauleap}

In this approximation method we build the matrices
\[
(B_{r})_{i,j}=\left\{ 
\begin{array}{rl}
-a_{r}(x_{j}) & \mbox{if $i=j$}, \\ 
a_{r}(x_{j}) & \mbox{if $i-j  =d_{r}$}, \\
0 & \textnormal{otherwise}
\end{array}
\right.
\]
where $a_{r}(x)$ are the propensity functions. The matrix $A$ in \eqref{eq_exact} can be written as
\[
A=\sum_{r=1}^M B_{r}\,.
\]
The solution of the linear CME \eqref{eq_exact2} can be approximated by
\begin{equation}  \label{eq_accelerated1}
\mathcal{P}\left(\bar{t}+\tau\right)=e^{\tau A}\cdot \mathcal{P}\left(\bar{t}\right) \approx   e^{\tau B_{1}}
e^{\tau B_{2}} \cdots e^{\tau B_{M}} \cdot P\left(\bar{t}\right)\,.
\end{equation}
Note that the evolution of state probability by $e^{\tau B_{j}}\cdot P\left(\bar{t}\right)$ describes the change in probability when only reaction $j$ fires in the time interval $\tau$. The corresponding evolution of the number of molecules that samples the evolved
probability is
\[
 x\left(\bar{t}+\tau\right)=x\left(\bar{t}\right)+V_{j}\, K\left(a_j\left(x\left(\bar{t}\right)\right) \tau\right).
\]
where $K\left(a_j\left(x\left(\bar{t}\right)\right) \tau\right)$ is a random number drawn from a Poisson distribution with parameter  $a_j\left(x\left(\bar{t}\right)\right) \tau$, and $V_{j}$ is the $j$-{th} column of stoichiometry matrix. 

The approximate solution \eqref{eq_accelerated1} accounts for the change in probability due to a sequential firing
of reactions $M$, $M-1$, down to $1$. Sampling from the resulting probability density can be done by 
changing the system state sequentially consistent with the firing of each reaction. This results in the following
accelerated tau-leaping algorithm:
\begin{equation}\label{eq_accelerated2}
\begin{array}{l}
\hat{X}_{M}= x\left(\bar{t}\right) \\
\textnormal{for  } i=M,M-1,\cdots,1 \\
\qquad \hat{X}_{i-1}=\hat{X}_{i}+V_{i}\, K\left(a_{i}\left(\hat{X}_{i}\right)\tau\right) \\
x(\bar{t}+\tau)=\hat{X}_{0}.  
\end{array}
\end{equation}

Moreover, \eqref{eq_accelerated1} can also be written as:
\begin{eqnarray}  \label{eq_accelerated3} 
\mathcal{P}\left(\bar{t}+\tau\right) \approx  e^{\tau B_{1}}
e^{\tau B_{2}} \cdots e^{\tau B_{M}} \cdot P\left(\bar{t}\right)\,\\ \nonumber
 \approx \left( e^{\tau B_{1}} 
e^{\tau B_{2}} \cdots e^{\tau B_{\frac{M}{2}-1}} \right) \cdot \\  \nonumber
 \left(e^{\tau B_{\frac{M}{2}}} e^{\tau B_{\frac{M}{2}+1}} \cdots e^{\tau B_{M}} \cdot P\left(\bar{t}\right) \right).\ \nonumber
\end{eqnarray}
Then, \eqref{eq_accelerated2} can be written as:
\begin{equation} \label{eq_accelerated4} 
\begin{array}{l}
\hat{X}_{M}= x\left(\bar{t}\right) \\
\textnormal{for  } i=M,M-1,\cdots,\frac{M}{2} \\
\qquad \hat{X}_{i-1}=\hat{X}_{i}+V_{i}\, K\left(a\left(\hat{X}_{M}\right)\tau\right) \\
\textnormal{for  } i=\frac{M}{2}-1,\cdots,1 \\
\qquad \hat{X}_{i-1}=\hat{X}_{i}+V_{i}\, K\left(a\left(\hat{X}_{\frac{M}{2}-1}\right)\tau\right) \\
x(\bar{t}+\tau)=\hat{X}_{0}.  
\end{array}
\end{equation}
%

\subsection{Symmetric accelerated tau-leaping}
\label{sect:SymmetricTauleap}
A more accurate version of accelerated tau-leaping can be constructed by using symmetric Strang
splitting \eqref{eq_strang} to approximate the matrix exponential in \eqref{eq_accelerated1}.
Following the procedure used to derive \eqref{eq_accelerated2} leads to the following sampling algorithm:
\begin{equation}\label{eq_symmetric}
\begin{array}{l}
\hat{X}_{M}= x\left(\bar{t}\right) \\
\textnormal{for  } i=M,M-1,\cdots,1 \\
\qquad \hat{X}_{i-1}=\hat{X}_{i}+V_{i}\, K\left(a_{i}\left(\hat{X}_{i}\right)\tau/2\right) \\
\textnormal{for  } i=1,2,\cdots,M \\
\qquad \hat{X}_{i}=\hat{X}_{i}+V_{i-1}\, K\left(a_{i}\left(\hat{X}_{i-1}\right)\tau/2\right) \\
x(\bar{t}+\tau)=\hat{X}_{M}.  
\end{array}
\end{equation}
%
\section{Numerical experiments}
\label{sect:NumericalExperim}
The above approximation techniques are used to solve two test systems,  reversible isomer and the Schlogl reactions \cite{Cao_2004_stability}. The experimental results are presented in following sections.                                                                                                                                                                                                                                                                                                                                                                                                                                                                                                                                                                                                                                                                                                                                                                                                                                                                                                                                                                                                                                                                                                                                                                                                                                                                                                                                                                                                                                                                                                                                                                                                                                                                                                                                                                                                                                                                                                                                                                                                                                                                                                                                                                                                                                                                                                                                                                                                                                                                                                                                                                                                                                                                                                                                                                                                
\subsection{Isomer reaction}
\label{sect:Isomer}

The reversible isomer reaction system is \cite{Cao_2004_stability}
\begin{equation}
\label{eqn:isomer}
\ce{
x_1 <=>[\ce{c_1}][\ce{c_2}]  x_2.
}
\end{equation}
The stoichiometry matrix and the propensity functions are: 
\[
V= 
\left[\begin{array}{rr}
-1 & 1 \\
1 & -1 
\end{array}\right]\,, \qquad
\begin{array}{l}
a_{1}(x)= c_{1}x_{1} \,, ~~\\
a_{2}(x) = c_{2}x_{2} \,.
\end{array}
\]
The reaction rate values are $ c_{1}=10$, $c_{2}=10$ (units), the time interval is $[0,T]$ with $T=10$ (time units), initial conditions are $x_{1}(0)=40$,
$x_{2}(0)=40$ molecules, and maximum values of species are $Q^1=80$ and $Q^2=80$ molecules.

\begin{figure}[tb]
	\begin{centering}
	\subfigure[10,000  SSA runs versus the exact solution \eqref{eq_exact2}]{
	\includegraphics[width=0.47\textwidth,height=0.3\textwidth]{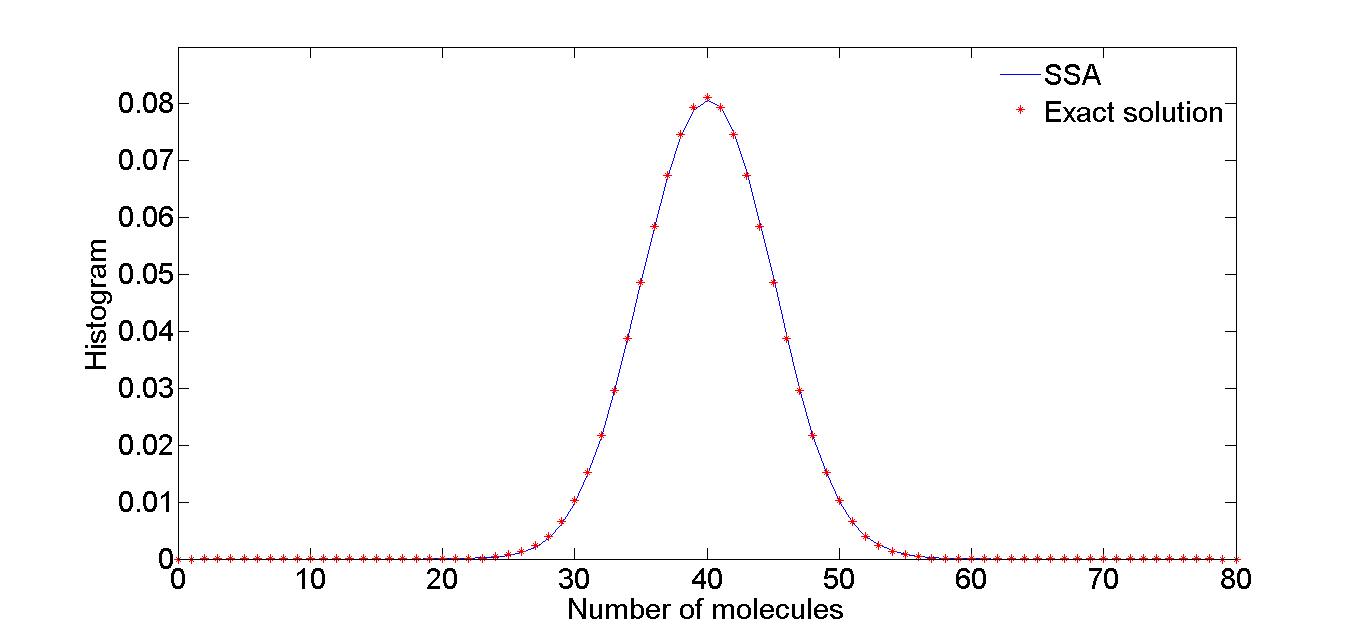} 
	\label{fig:exact_isomer}
	}
	\subfigure[Exact solution \eqref{eq_exact2} versus the approximation to exact solution using sum of exponentials \eqref{eq_sumofexponent3} ]{	
	\includegraphics[width=0.47\textwidth,height=0.3\textwidth]{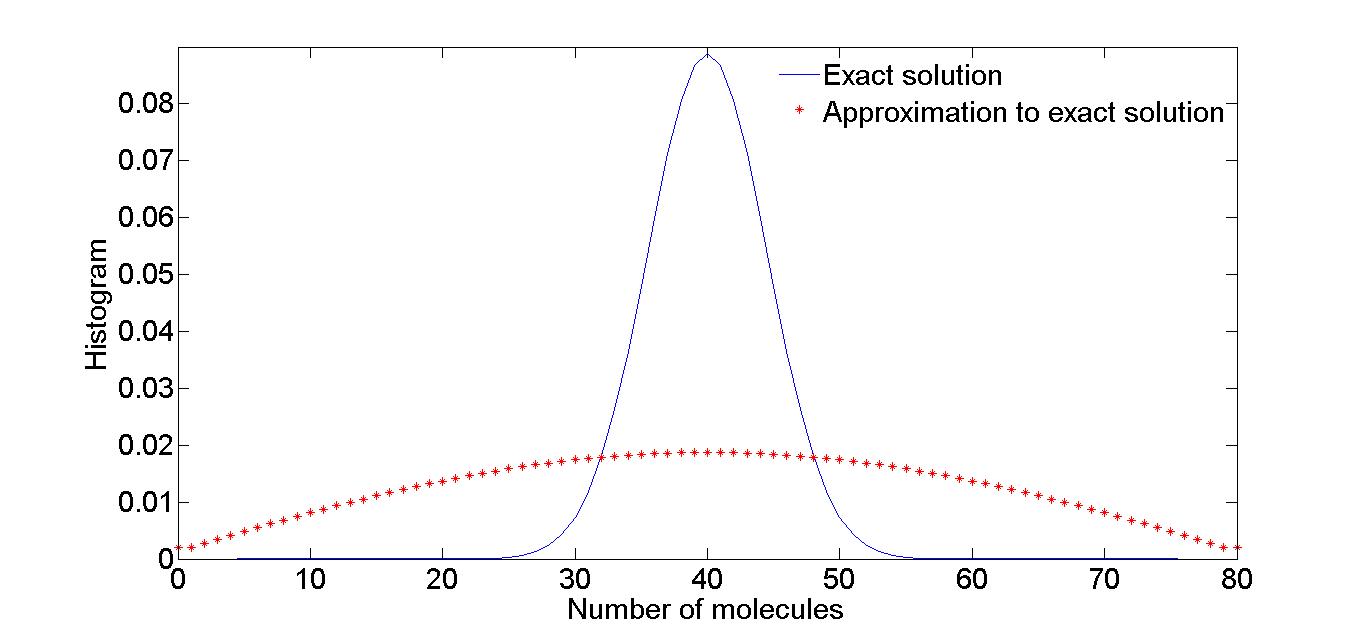} 
	\label{fig:approx_isomer}
        }
	\caption{Histograms of the isomer system \eqref{eqn:isomer} results at the final time T=10.}
	\label{fig:isomer}
	\end{centering}
\end{figure}

The exact exponential solution of CME obtained from \eqref{eq_exact2} is a joint probability  distribution vector for 
the two species at final time. Figure \ref{fig:exact_isomer} shows that the histogram of 10,000 SSA solutions is very close to the exact exponential solution. The approximate solution using the sum of exponentials \eqref{eq_sumofexponent3} is illustrated in Figure \ref{fig:approx_isomer}. This approximation is not very accurate since it uses only the current state of the system. Other approximation methods based on the product of exponentials \eqref{eq__productofexponents} and Strang splitting \eqref{eq_strang} are not very strong approximations as the exact solution hence, the results are not reported.

\begin{figure}[tb]
	\begin{centering}
	\includegraphics[width=0.475\textwidth,height=0.3\textwidth]{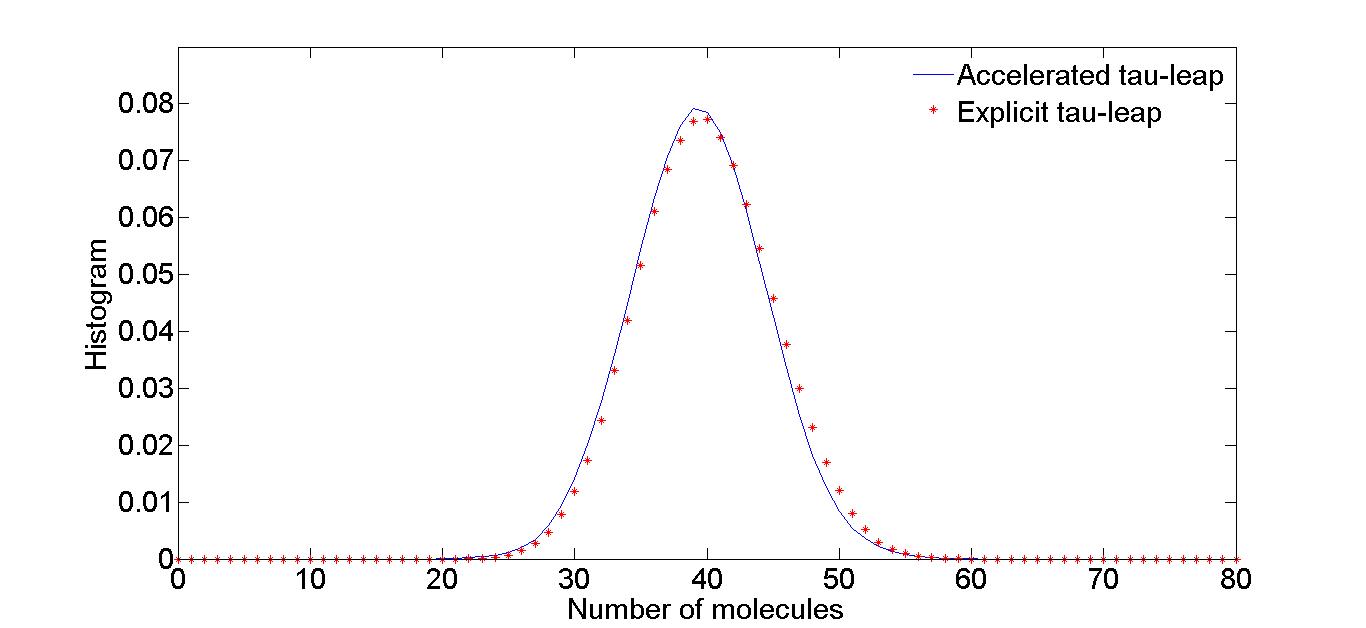} 
	\caption{Isomer system \eqref{eqn:isomer} solutions provided by the traditional tau-leap \eqref{eq_tauleap} and by accelerated tau-leap \eqref{eq_accelerated2} methods at the final time T=10 (units). A small time step of $\tau=0.01$ (units) is used. The number of samples for both methods is 10,000.}
	\label{fig:accelerated_isomer}
	\end{centering}
\end{figure}

\begin{figure}[tb]
	\begin{centering}
	\subfigure[$\tau = 0.01$ (units)]{
	\includegraphics[width=0.475\textwidth,height=0.3\textwidth]{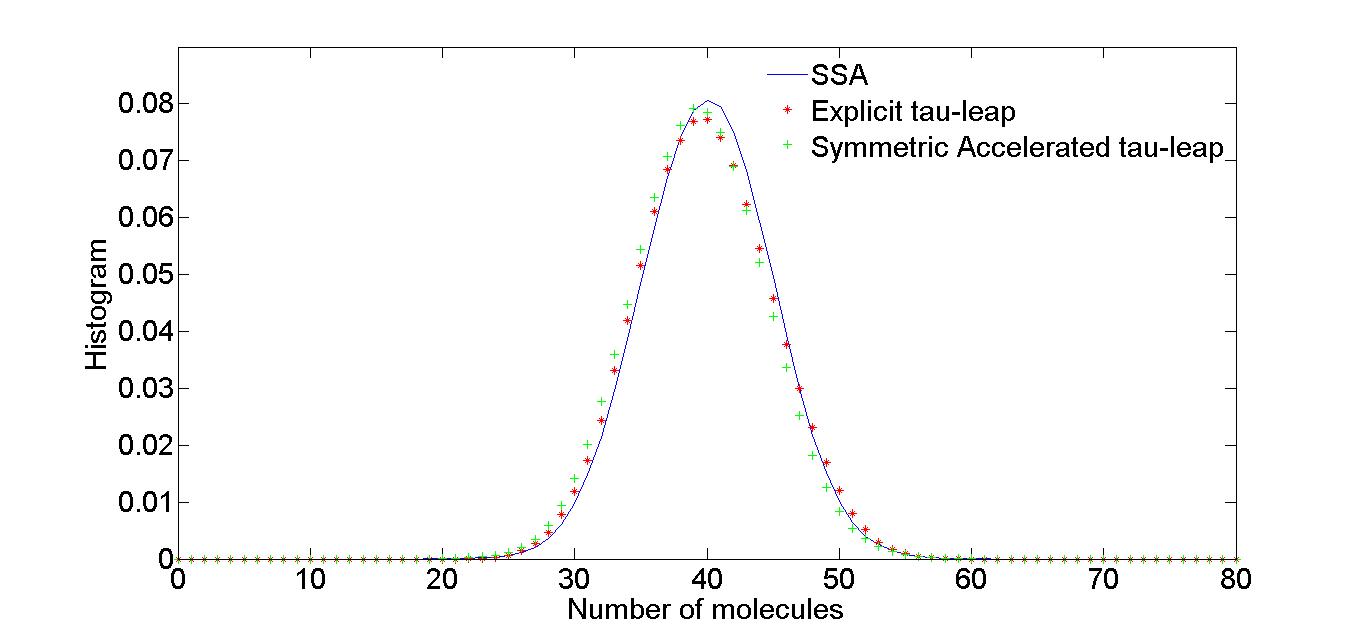} 
	\label{fig:symmetric_isomer_1per}
	}
	\subfigure[$\tau = 0.1$ (units)]{	
	\includegraphics[width=0.475\textwidth,height=0.3\textwidth]{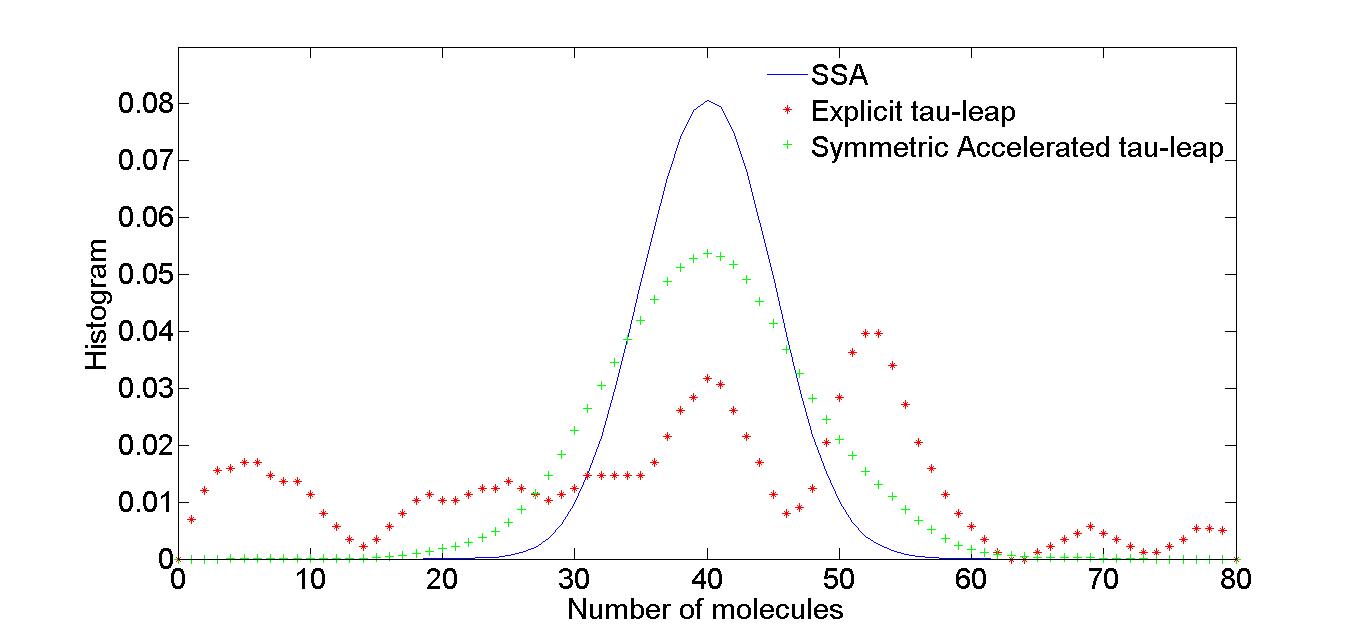} 
	\label{fig:symmetric_isomer_1en}
        }
	\caption{Histograms of isomer system \eqref{eqn:isomer} solutions obtained with SSA, traditional tau-leap  \eqref{eq_tauleap}, and symmetric accelerated tau-leap \eqref{eq_symmetric}
methods at the final time T=10. The number of samples is 10,000  for all methods.}
	\label{fig:symmetric_isomer_1ten}
	\end{centering}
\end{figure}

The results  reported in Figure \ref{fig:accelerated_isomer} indicate that for small time steps $\tau$ the accelerated tau-leap \eqref{eq_accelerated2} solution is very close to the results provided by traditional explicit tau-leap. Symmetric accelerated tau-leap method \eqref{eq_symmetric} yields even better results, as shown in Figure \ref{fig:symmetric_isomer_1ten}.
For small time steps the traditional and symmetric accelerated methods give similar results,
however, for large time steps, the results of the  symmetric accelerated method is considerably more stable.

%
%

\subsection{Schlogl reaction}
\label{sect:Schlogl}

We next consider the Schlogl reaction system \cite{Cao_2004_stability}

\begin{equation}
\label{eqn:schlogl}
\begin{array}{r}
\ce{
B_{1} + 2x  <=>[\ce{c_1}][\ce{c_2}]  3x
}\\
\ce{
B_{2} <=>[\ce{c_3}][\ce{c_4}]  x
}
\end{array}
\end{equation}
whose solution has a bi-stable distribution. Let $N_1$, $N_2$ be the numbers of molecules of species $B_1$ and $B_2$, respectively.
The reaction stoichiometry matrix and the propensity functions are: 
\[
\begin{array}{l}
V= 
\begin{bmatrix}
1 & -1 & 1 & -1
\end{bmatrix} \\
\begin{array}{l}
a_{1}(x)= \frac{c_{1}}{2}N_{1}x(x-1), \\ 
a_{2}(x) = \frac{c_{2}}{6}N_{1}x(x-1)(x-2), \\ 
a_{3}(x) = c_{3}N_{2}, \\ 
a_{4}(x) = c_{4}x. 
\end{array}%
\end{array}
\]
The following parameter values (each in appropriate units) are used:
\begin{small}
\[
\begin{array}{lll}
c_{1}=3 \times 10^{-7}, &c_{2}=10^{-4}, &c_{3}=10^{-3}, \\
c_{4}=3.5, &N_{1}=1 \times 10^5, &N_{2}=2 \times10^5.
\end{array}%
\]
\end{small}
with the final time $T=4$ (units), the initial condition $x(0)=250$ molecules, and the maximum values of species $Q^1=900$ molecules.

\begin{figure}[tb]
	\begin{centering}
	\subfigure[10,000  SSA runs versus the exact solution \eqref{eq_exact2}]{
	\includegraphics[width=0.475\textwidth,height=0.3\textwidth]{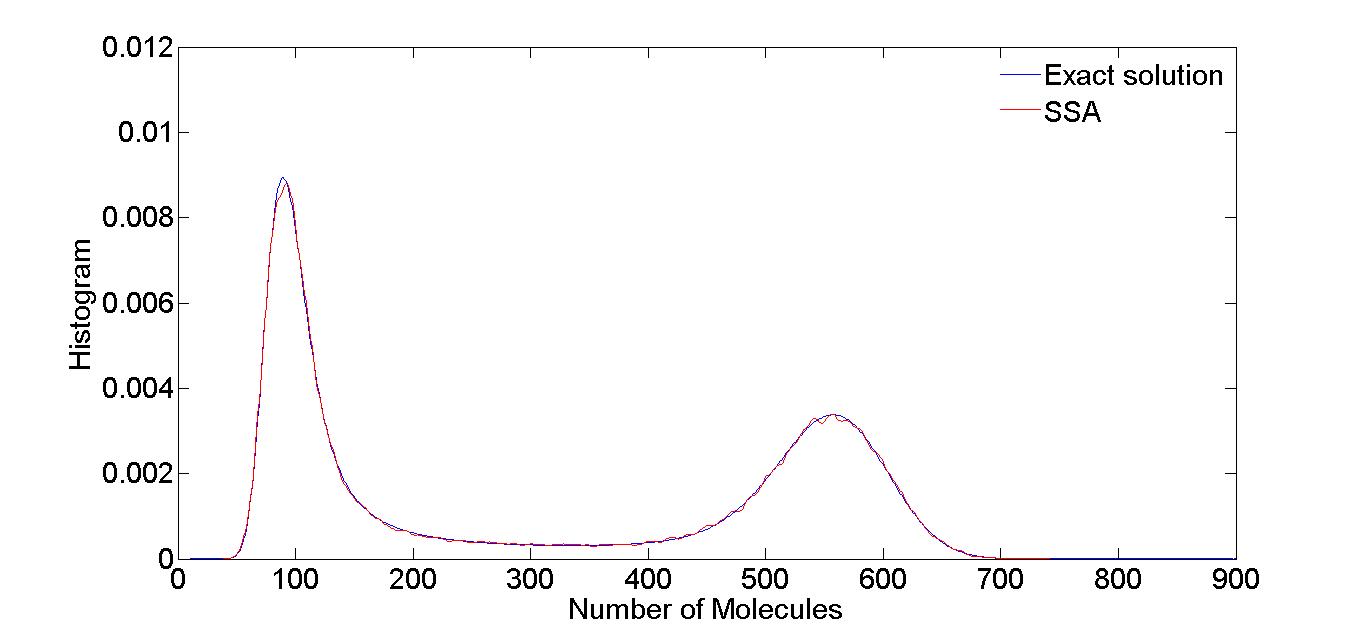} 
	\label{fig:exact_schlogl}
	}
	\subfigure[Exact solution \eqref{eq_exact2} versus the approximation to exact solution using sum of exponentials \eqref{eq_sumofexponent3}]{	
	\includegraphics[width=0.475\textwidth,height=0.3\textwidth]{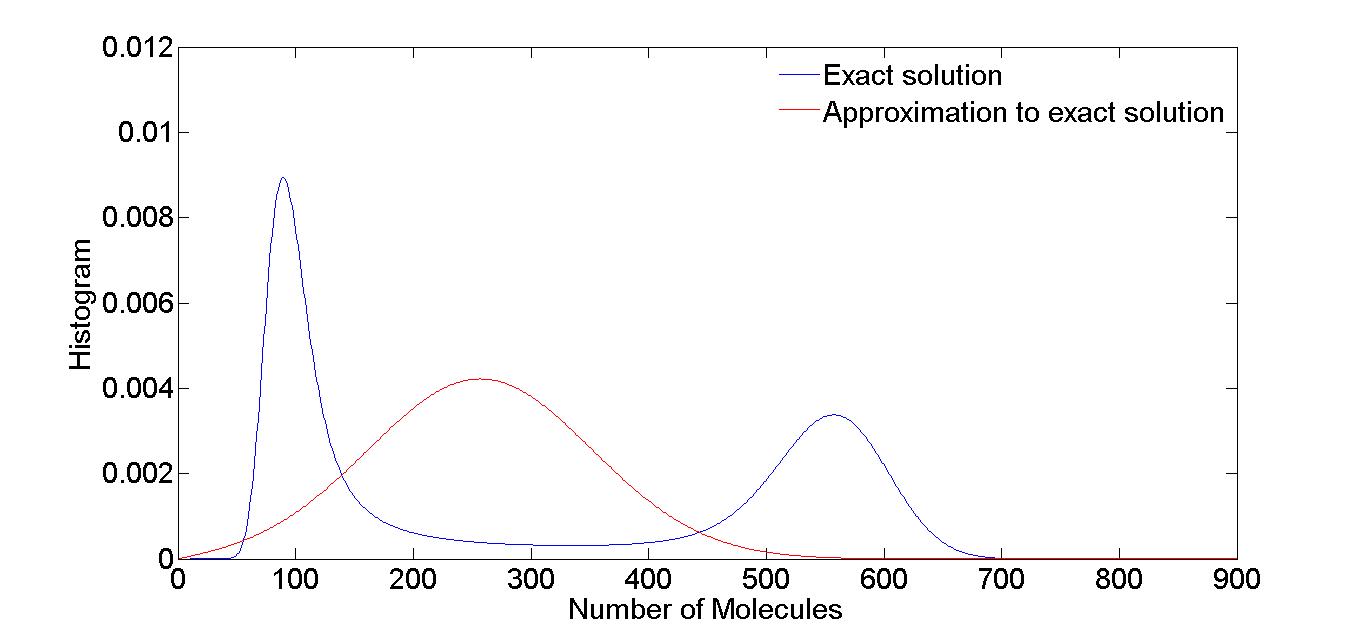} 
\label{fig:approx_exact_schlogl}}
	\caption{Histograms of Schlogl system \eqref{eqn:schlogl} results at final time T=4 (units).}
	\end{centering}
\end{figure}

Figure \ref{fig:exact_schlogl} illustrates the result of exact exponential solution \eqref{eq_exact2} versus SSA. Figure \ref{fig:approx_exact_schlogl} reports the sum of  exponentials \eqref{eq_sumofexponent3} result which is not a very good approximation. The product of exponentials \eqref{eq__productofexponents} and Strang splitting \eqref{eq_strang} results are not reported here since they are poor in approximation. 

\begin{figure}[tb]
	\begin{centering}
	\subfigure[Traditional tau-leap  \eqref{eq_tauleap} and accelerated tau-leap \eqref{eq_accelerated2}]{
	\includegraphics[width=0.475\textwidth,height=0.3\textwidth]{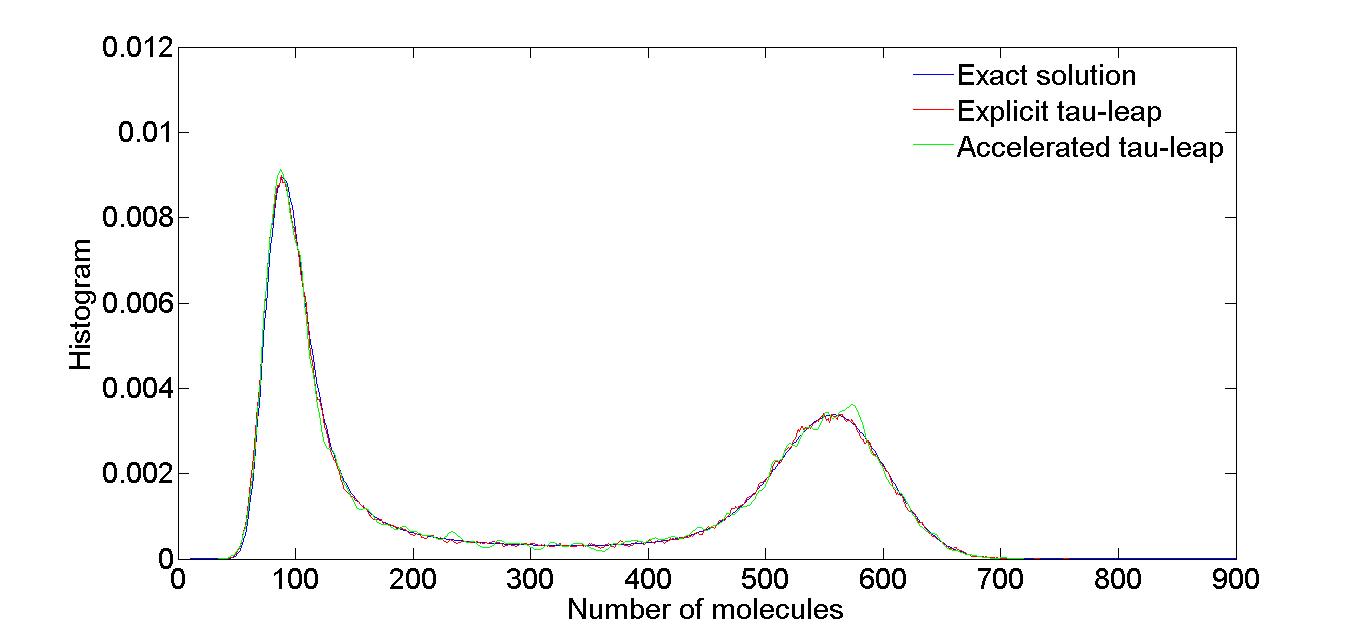} 
	\label{fig:accelerated_schlogl}
	}
	\subfigure[Traditional tau-leap  \eqref{eq_tauleap} and symmetric accelerated tau-leap \eqref{eq_symmetric}]{	
	\includegraphics[width=0.475\textwidth,height=0.3\textwidth]{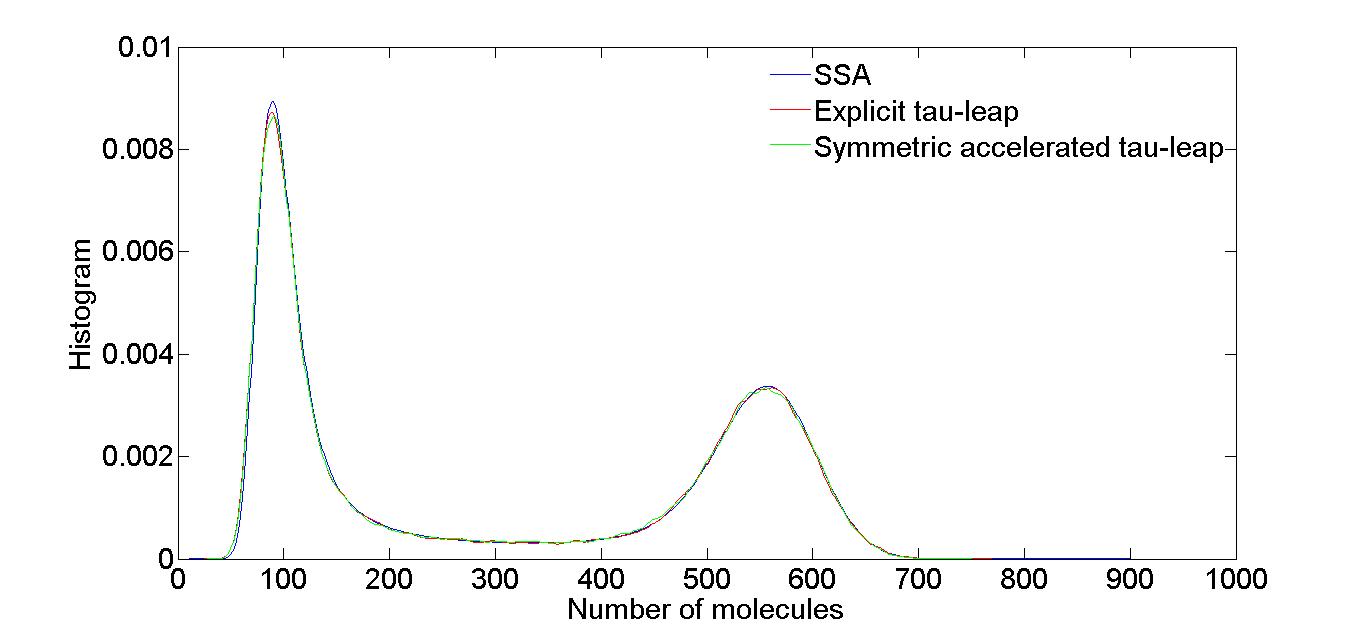} 
\label{fig:symmetric_schlogl}}
	\caption{Histograms of Schlogl system \eqref{eqn:schlogl} solutions with $\tau=0.0001$ (units), final time T=4 (units), and 10,000 samples.}
	\end{centering}
\end{figure}

Figures \ref{fig:accelerated_schlogl} and \ref{fig:symmetric_schlogl} present the results obtained with the accelerated tau-leap and the symmetric tau-leap, respectively. For small time step the results are very accurate. However, for large step sizes, the results quickly become less accurate. The lower accuracy may affect systems having more reactions. The accuracy can be improved to some extent using the strategies described in \eqref{eq_accelerated3} and \eqref{eq_accelerated4}. 

%
%
\section{Conclusions}
\label{sect:Conclusion}
This study proposes new numerical solvers for stochastic simulations of chemical kinetics.
The proposed approach exploits the linearity of the CME and the exponential form of its exact solution. 
The matrix exponential appearing in the CME solution is approximated as a product of simpler matrix exponentials.
This leads to an approximate (``numerical'') solution of the probability density evolved to a future time.
The solution algorithms sample exactly this approximate probability density and provide extensions of the traditional tau-leap approach.

Different approximations of the matrix exponential lead to different numerical algorithms: Strang splitting, column splitting, accelerated tau-leap, and symmetric accelerated tau-leap. Current work by the authors focuses on improving the accuracy of these novel approximation techniques for stochastic chemical kinetics.
\bibliographystyle{plain}
\bibliography{main.bib}

%
\ifx

\fi
\appendix
\section{Example}
\label{sec:example}
We exemplify the process of building matrix A \eqref{eq_exact} for the Schlogl and isomer reactions. 

\subsection{Isomer reaction}
Here for simplicity, we exemplify the implementation of the system for the maximum values of species $Q^1 = 2$ and $Q^2 = 2$. According to \eqref{eq_Q}, $Q=(Q^{1}+1) \times (Q^{2}+1)=3^2$.

The vector $d$ according to \eqref{eq_d} is $[2,-2]$. The state
matrix which contains all possible states has dimension $81^2 \times 2$ matrix:
\[
\mathbf{x} = 
\begin{bmatrix}
\ 0 & 1 &2 & 0 & 1 &2 & 0 & 1 &2  \\[0.3em] 
\ 0 & 0 &0 & 1 & 1 &1 & 2 & 2 &2 
\end{bmatrix}^\top \in \mathbb{R}^{3^2 \times 2}.
\]
The matrix $\mathbf{A} \in \mathbb{R}^{Q\cdot Q \times Q \cdot Q}$
As an example for a maximum number of species $Q^{1}=2$, $Q^{2}=2$
the matrix $\mathbf{A}$ is: 

\begin{center}
\begin{small}
\begin{eqnarray} \nonumber
\mathbf{A} = 
\begin{bmatrix}
-a_{0}(\mathbf{x}_{1,:}) & 0 & a_{2}(\mathbf{x}_{3,:}) & 0 & 0 & 0 & 0 & 0 & 0 \\[0.3em] 
0 & -a_{0}(\mathbf{x}_{2,:}) & 0 & \ddots & 0 & 0 & 0 & 0 & 0 \\[0.3em] 
a_{1}(\mathbf{x}_{1,:}) & 0 & -a_{0}(\mathbf{x}_{3,:}) & 0 & \ddots & 0 & 0 & 0 & 0 \\[0.3em] 
0 & a_{1}(\mathbf{x}_{2,:}) & 0 & \ddots & 0 & \ddots & 0 & 0 & 0 \\[0.3em] 
0 & 0 & a_{1}(\mathbf{x}_{3,:}) & 0 & \ddots & 0 & a_{2}(\mathbf{x}_{7,:}) & 0 & 0 \\[0.3em] 
0 & 0 & 0 & \ddots & 0 & \ddots & 0 & a_{2}(\mathbf{x}_{8,:}) & 0 \\[0.3em] 
0 & 0 & 0 & 0 & \ddots & 0 & -a_{0}(\mathbf{x}_{7,:}) & 0 & a_{2}(\mathbf{x}_{9,:}) \\[0.3em] 
0 & 0 & 0 & 0 & 0 & \ddots & 0 & -a_{0}(\mathbf{x}_{8,:}) & 0 \\[0.3em] 
0 & 0 & 0 & 0 & 0 & 0 & a_{1}(\mathbf{x}_{7,:}) & 0 & -a_{0}(\mathbf{x}_{9,:}) 
\end{bmatrix}  \in \mathbb{R}^{9\times 9}\,.
\end{eqnarray}
\end{small}
\end{center}

\subsection{Schlogl reaction}
Here for simplicity, we exemplify the implementation of the system for the maximum value of the number of molecules $Q^1=5$. 
According to 
\eqref{eq_Q} the dimensions of A are: $\left(Q^1+1 \times Q^1+1\right)= 6 \times 6$. The vector $d$ \eqref{eq_d} for this system $%
[1,-1,1,-1]$. All possible states for this system are contained in the state vector
\[
\mathbf{x} = [0,1,2, \cdots, 5]^\top  \in \mathbb{R}^{1 \times 6}.
\]
As an example matrix A for maximum number of molecules $Q=5$ is the following tridiagonal matrix:
\begin{center}
\begin{small}
\begin{eqnarray} \nonumber
\mathbf{A} = 
\begin{bmatrix} 
-a_{0}(\mathbf{x}_1) & a_{2}(\mathbf{x}_2)+a_{4}(\mathbf{x}_2) & 0 & 0 & 0 & 0 \\[0.3em] 
a_{1}(\mathbf{x}_1)+a_{3}(\mathbf{x}_1) & -a_{0}(\mathbf{x}_2) & \ddots & 0 & 0 & 0 \\[0.3em] 
0 & a_{1}(\mathbf{x}_2) + a_{3}(\mathbf{x}_2) & \ddots & \ddots & 0 & 0 \\[0.3em] 
0 & 0 & \ddots & \ddots & a_{2}(\mathbf{x}_5) + a_{4}(\mathbf{x}_5) & 0 \\[0.3em] 
0 & 0 & 0 & \ddots & -a_{0}(\mathbf{x}_5) & a_{2}(\mathbf{x}_6) + a_{4}(\mathbf{x}_6) \\[0.3em] 
0 & 0 & 0 & 0 & a_{1}(\mathbf{x}_5) + a_{3}(\mathbf{x}_5) & -a_{0}(\mathbf{x}_6) 
\end{bmatrix}  \in \mathbb{R}^{6 \times 6}.
\end{eqnarray}
\end{small}
\end{center}

\label{lastpage}
\end{document}